\documentclass[10pt]{article}

\usepackage[OT2,T1]{fontenc}
\usepackage{amssymb,amsfonts,amsmath,amsthm,amscd}
\usepackage{relsize}
\usepackage{authblk}
\usepackage[margin=2.8cm]{geometry}
\usepackage{bigstrut}
\usepackage{mathtools}
\usepackage{titlefoot}
\usepackage{cite}
\usepackage{xcolor}
\usepackage{subcaption}
\usepackage[inline]{enumitem}
\usepackage{float}
\usepackage{subfig}

\DeclareMathOperator\sgn{sgn}
\DeclareMathOperator{\ctan}{ctan}

\DeclareMathOperator{\am}{am}

\DeclareMathOperator{\dn}{dn}
\DeclareMathOperator{\nd}{nd}
\DeclareMathOperator{\sn}{sn}
\DeclareMathOperator{\cs}{cs}
\DeclareMathOperator{\sc1}{sc}

\DeclareMathOperator{\nc}{nc}
\DeclareMathOperator{\cn}{cn}
\DeclareMathOperator{\ns}{ns}
\DeclareMathOperator{\cd}{cd}

\def\RR{\mathbb{R}}

\def\t{\mathsmaller{T}}

\def\H{H_{5}}
\def\h{\mathfrak{h}_5}
\def\DD{\mathcal{D}}

\newtheorem{lemma}{{Lemma}}[section]
\newtheorem{theorem}{{Theorem}}[section]
\newtheorem{proposition}{{Proposition}}[section]
\newtheorem{cor}{{Corollary}}[section]
\theoremstyle{definition}
\newtheorem{remark}{{Remark}}[section]

\newtheorem{examplex}{Example}[section]
\newenvironment{example}
  {\pushQED{\qed}\examplex}
  {\popQED\endexamplex}

\title{Classification of the
geodesic curves of a sub-Riemannian LR system on the Heisenberg group in dimension 5}
\author{Milan Pavlovi\'{c}}
\affil{\small{University of Belgrade, Faculty of Mathematics, Belgrade, Serbia}}
\date{\today}

\begin{document}
\maketitle\unmarkedfntext{2020 \emph{Mathematics Subject Classification}. 53C17, 22E25, 37J35\\ 
\indent\phantom{K}\emph{Key words and phrases}: Heisenberg group, left-invariant metrics, sub-Riemannian metric, right-invariant distributions, geodesic flow, completely integrable Hamiltonian system.}

\begin{abstract}
We study sub-Riemannian structures obtained by the restriction of the
standard left-invariant Riemannian metric to a right-invariant completely nonholonomic
distribution on the second Heisenberg group.
The corresponding Hamiltonian system is completely integrable and in this paper we study its solutions. We obtain a complete classification of the geodesic curves by studying the reduced Hamiltonian system. Moreover, we compute $r(t)$ as the first step of the quadrature in the reduced system in the hyper-spherical coordinates and present numerical examples that illustrate the theoretical results.
\end{abstract}

\section*{Introduction}
Sub-Riemannian geometry can be described as geometry with restrictions. Motion is only possible in certain directions (which vary from point to point). These constraints are defined by the distribution in which the length is measured. If the distribution has a certain 'nice' property, every point in space is reachable. There are several recent papers presenting new integrable sub-Riemannian geodesic flows on Lie
groups and homogeneous spaces, such as~\cite{
BDAB2025, BJ2025, JSV2025, VJ2025, PS2024}. In this paper we will study one such integrable sub-Riemannian geodesic flow, found in~\cite{PS2024}.  In general, sub-Riemannian geometry on Lie groups is studied by considering metric and the distribution, both of which are left-invariant. Less studied are situations where the metric is left-invariant and the distribution is right-invariant, the so-called LR system. Such a situation leads to strongly non-linear dynamics. LR systems have been already studied in nonholonomic mechanics. The basic example is a
motion of a rigid body around a fixed point in the presence of a nonholonomic constraint: the
projection of the angular velocity of the body onto a fixed axis in the space is equal to zero (~\cite{FJ2004,VV1988}). LR systems are also studied sub-Riemannian framework (~\cite{Taim1997, AB2017, Maz2014}). In \cite{Taim1997} a sub-Riemannian LR system on a three-dimensional Heisenberg group is studied and it is found that such a system is both integrable and associated with a certain magnetic field, while in \cite{AB2017} geodesics of a sub-Riemannian LR system on the first Heisenberg group are classified. In \cite{PS2024} first results on sub-Riemannian LR systems in dimensions larger than three are given. In particular, the integrability of sub-Riemannian LR systems in the general Heisenberg group $H_{2n+1}$ is proved and the equivalence to motion in the presence of a higher-dimensional magnetic field is found. If $n > 2$, one has integrability in the non-commutative sense, while for $n=1$ and $n=2$, i.e. in dimensions three and five, systems are integrable in the classical, commutative sense. 

The integrability of sub-Riemannian LR systems on the Heisenberg group has been the subject of previous investigations. Motivated by these results, we turn our attention to the study of geodesic lines, whose analysis provides a natural continuation of this research and contributes to a more comprehensive understanding of the underlying sub-Riemannian geometry. We classify the sub-Riemannian geodesics in dimension $5$ and we study the reduced Hamiltonian system in hyperspherical coordinates. This coordinates allows us to classify geodesics according to the properties of $r(t)$. 

The paper is organized as follows. Section~\ref{sec1} reviews the necessary background on sub-Riemannian geometry, the five-dimensional Heisenberg group, and elliptic integrals, while also introducing the notation used throughout the paper. In Section~\ref{sec2}, we revisit the integrability of a sub-Riemannian LR system on $\H$, previously established in~\cite{PS2024} for Heisenberg groups of arbitrary dimension. We derive the corresponding reduced Hamiltonian system and present its first integrals. In Section~\ref{sec3} a classification of the trajectories of the reduced Hamiltonian system in hyper-spherical coordinates is obtained. This classification extends naturally to the classification of the geodesic curves in the second Heisenberg group, and is the main result of the paper. We conclude this section by deriving explicit formulae for the first step of the quadrature integration, $r=r(t)$, and by presenting numerical examples that illustrate the theoretical results.

\section{Preliminaries}\label{sec1}
In this section, we introduce the notation as well as some well-known facts of sub-Riemannian geometry, most of which can be found for example in~\cite{ABB2019, Mont2002,Str1986, Taim1997, VJ2016}, Heisenberg group, and Elliptic integrals.

\subsection{Sub-Riemannian geometry}

Let $Q$ be a manifold of dimension $n$, $\mathcal{D}$ a distribution of constant rank $m$ of the tangent bundle $TQ$ and $\langle \cdot,\cdot \rangle_{sR}(q)$ a smoothly varying positive definite bilinear form acting as an inner product on this distribution. A sub-Riemannian structure on $Q$ is a pair $(\mathcal{D}, \langle \cdot,\cdot \rangle_{sR})$. If at each point $q\in Q$ vector fields of $\mathcal{D}_q$ together with all their Lie brackets span the entire $T_qQ$, we say that this distribution is completely nonholonomic (or bracket-generating).

A horizontal (or admissible) curve is a curve whose tangent vector lies in the distribution. More precisely, the absolutely continuous curve $\gamma: [0,1]\rightarrow Q$ is horizontal if $\dot{\gamma}(t)\in \mathcal{D}_{\gamma(t)}$ for almost every $t\in[0,1]$. For horizontal $\gamma$ we define its length as
\begin{align*}
 \ell(\gamma)=\int_0^1 \|\dot{\gamma}(t)\|dt,
\end{align*}
where $\|\dot{\gamma}\|=\sqrt{\langle \dot{\gamma}, \dot{\gamma}\rangle_{sR}}$. In this case, the sub-Riemannian or Carnot-Carath\'{e}odory distance $dist_{CC}(x,y)=\inf \ \ell (\gamma)$ is defined, and the infimum is taken over all admissible curves $\gamma$ connecting $x$ and $y$. If there is no such curve, the Carnot-Carath\'{e}odory distance is infinite. For distributions with bracket-generating property, the theorems of Chow~\cite{Chow1939} and Rashievskii~\cite{Ras1938} guarantee that there is always an admissible curve connecting any two points.

For a given $\langle \cdot,\cdot \rangle_{sR}(q)$, the linear mapping $g(q):T_q^*M\rightarrow \mathcal{D}_q \subset T_q M$ can be defined as follows: $g(q)\lambda$ is a unique vector such that
\begin{equation*}
 \langle g(q)\lambda,Y \rangle_{sR}(q)=\lambda(Y(q)) ,\quad \forall Y\in\mathcal{D}_q,
\end{equation*}
with the usual pairing of vectors and covectors right-hand side. The symmetric tensor $g^{ij}$ is called Carnot-Carathdodory metric tensor. It is non-negative definite, but generally not positive definite, since it is not onto.

For the pairing $\lambda(X(q))$, $q\in Q$ between the fixed vector field $X$ and the covector $\lambda$ we write $P_X(q,\lambda)$. The function $P_X\colon T^*M\to\mathbb R$  is called the momentum function. Let $X_1,\dots,X_m$ be a local basis
of vector fields of the distribution, in coordinates $X_a=\sum_{j} X_a^j(q)\frac{\partial}{\partial x_j}$. The momentum functions are given by
{$P_{X_a}=\sum_{j} X_a^j(q)\lambda_j$}. Let $g_{ab}(q)$ be the matrix of the inner products of the vectors $X_a$ and $g^{ab}(q)$ their inverse matrix. Then $g^{ab}(q)$ is an $m\times m$ matrix-valued function (where $m$ is the rank of the distribution) defined in an open set of $Q$. The sub-Riemannian Hamiltonian function can be written as
\begin{align*}
 H(q,\lambda)=\frac{1}{2}\sum_{a,b} g^{ab}(q)P_{X_a}(q,\lambda)P_{X_b}(q,\lambda).
\end{align*}
With $\{f,g\}=X_g(f)$ we denote the Poisson bracket of $f,g\in T^*Q$, where $X_g$ is the Hamiltonian vector field associated with $g$. In canonical coordinates, the Poisson bracket is given by:
\begin{align*}
 \{f,g\}=\sum_i \frac{\partial f}{\partial q_i}\frac{\partial g}{\partial \lambda_i}- \frac{\partial f}{\partial \lambda_i}\frac{\partial g}{\partial q_i}.
\end{align*}
The corresponding Hamilton equations on the cotangent bundle are

\begin{align}\label{eq:HamEq}
 \dot{q}_i=\{q_i,H\}=\frac{\partial H}{\partial \lambda_i},\qquad\dot{\lambda}_i=\{\lambda_i,H\}=-\frac{\partial H}{\partial q_i}.
\end{align}
The projections of solutions of the Hamiltonian system~\eqref{eq:HamEq} are locally length-minimising curves, which are referred to as normal geodesics (see~\cite{Taim1997,Str1986}). In contrast to Riemannian geometry, in sub-Riemannian geometry, we have locally minimising curves that do not arise as solutions of the Hamiltonian system, abnormal geodesics. The horizontal distribution in the Heisenberg group is the contact distribution, and all geodesics are solutions of the Hamiltonian system~{\cite[Proposition 4.38]{ABB2019}.}

\subsection{Heisenberg group $\H$}

The second Heisenberg group $\H$ is a five-dimensional two-step nilpotent Lie group constructed on the base manifold $\RR ^{4}\oplus \RR$ by multiplication
\begin{align*}
(u, \zeta )\cdot (v, \chi):= (u + v ,\, \zeta + \chi + \omega (u,v)).
\end{align*}
\begin{align*}
\omega (u,v) = u^\t J v,\quad \text{where } J=\begin{pmatrix}
	0 & 0 & -1 & 0\\
 0 & 0 & 0 & -1\\
 1 & 0 & 0 & 0 \\
	0 & 1 & 0 & 0
\end{pmatrix} .
\end{align*}
This group can also be seen as the group of matrices of the form
\begin{equation*}
\H = \left\{
\begin{pmatrix}
1 & x_1 & x_2 & z \\
0 & 1 & 0 & y_1 \\
0 & 0 & 1 & y_2 \\
0 & 0 & 0 & 1
\end{pmatrix}: x_1, x_2, y_1, y_2, z \in \mathbb{R} \right\},
\end{equation*}
with respect to the matrix multiplication. The corresponding Lie algebra $\h$ is spanned by $\frac{\partial}{\partial x_1},\frac{\partial}{\partial x_2},\frac{\partial}{\partial y_1}, \frac{\partial}{\partial y_2}$ and $\frac{\partial}{\partial z}$.
The commutator relations that are not zero are
\begin{equation}\label{commutators:Heisenberg}
[\frac{\partial}{\partial x_1},\frac{\partial}{\partial y_1}]=[\frac{\partial}{\partial x_2},\frac{\partial}{\partial y_2}]=\frac{\partial}{\partial z}.
\end{equation}
This representation is more convenient for our study and is used in this paper.

\subsection{Elliptic integrals}
In the following, elliptic integrals will be useful, therefore we introduce the notation here, for more about elliptic integrals and elliptic functions we refer to extensive literature (e.g. \cite{WV1996}, \cite{BF1971}). With $F$ and $E$ we denote elliptic integrals of the first and second kind respectively:
\begin{align*}
 F(\phi,k)&=\int_0^\phi\frac{dv}{\sqrt{1-k^2\sin^2{v}}}=u=\int_0^y \frac{dt}{\sqrt{(1-t^2)(1-k^2t^2)}}=\sn^{-1}(y,k), \quad [y=\sin{\phi}],\\
\text{and}\\
 E(\phi,k)&=\int_0^\phi\sqrt{1-k^2\sin^2{v}}dv=\int_0^u\dn^2{u_1}du_1=\int_0^y\sqrt{\frac{1-k^2t^2}{1-t^2}}dt=E(\am{u},k),
 \end{align*}
where $\phi=\am{u}$ is the amplitude, the parameter $0\leq k\leq 1$ is the modulus and $k'=\sqrt{1-k^2}$ is the complementary modulus.
Twelve elliptic Jacobian functions are denoted by $pq(u,k)$, where $p,q\in{\{s,c,d,n\}}$, e.g. elliptic sine is $\sn{(u,k)}$, elliptic cosine is $\cn{(u,k)}$ and delta amplitude is $\dn{(u,k)}$.

The exponentiation of the Jacobian function with $-1$ should always be understood as inversion. The abbreviation convention $E(u)$ should be understood as $E(u)=E(\am{u},k)$ if $k$ is clear.

\section{Integrability in $\H$}\label{sec2}

The group $\H$ acts on itself by left translations $L_g$: $L_g(h)=gh$ and by right translations $R_g$: $R_g(h)=hg$. Denote by $\DD_0$ the linear space spanned by $\{\frac{\partial}{\partial x_1},\frac{\partial}{\partial x_2},\frac{\partial}{\partial y_1}, \frac{\partial}{\partial y_2}\}$ at the identity. The left-invariant distribution generated by $\DD_0$ consists of $4$-planes $\DD_L=L_{g*}\DD_0$, where:
\begin{align*}
L_{g*}(\frac{\partial}{\partial x_k})=\frac{\partial}{\partial x_k},\quad L_{g*}(\frac{\partial}{\partial y_k})=\frac{\partial}{\partial y_k} + x_k \frac{\partial}{\partial z}, \quad k\in\{1,2\}\quad\text{and}\quad L_{g*}(\frac{\partial}{\partial z})=\frac{\partial}{\partial z},
\end{align*}
while the right-invariant distribution generated by $\DD_0$ is $\DD_R=R_{g*}\DD_0$, where:
\begin{align*}
R_{g*}(\frac{\partial}{\partial x_k})=\frac{\partial}{\partial x_k}+ y_k \frac{\partial}{\partial z},\quad R_{g*}(\frac{\partial}{\partial y_k})=\frac{\partial}{\partial y_k} , \quad k\in\{1,2\}\quad\text{and}\quad R_{g*}(\frac{\partial}{\partial z})=\frac{\partial}{\partial z}.
\end{align*}
We note that $\H$ is a contact manifold and $\DD_R$ is a contact distribution. From the relation~\eqref{commutators:Heisenberg} it follows that this distribution is bracket-generating.

In \cite{PS2024} the geodesic flow of the sub-Riemannian metric, which corresponds to the left-invariant metric in $H_5$ and the right-invariant distribution $\DD_R$, is considered and its integrability is proved. The left-invariant standard Riemannian metric in the coordinates $(x_1,x_2,y_1,y_2,z)$ is:
\begin{equation*}
   g_{ij}= \begin{pmatrix}
 1 & 0 & 0 & 0 & 0 \\
 0 & 1 & 0 & 0 & 0 \\
 0 & 0 & x_1^2+1 & x_1x_2 & -x_1 \\
 0 & 0 & x_1x_2 & x_2^2+1 & -x_2 \\
 0 & 0 & -x_1 & -x_2 & 1 
\end{pmatrix}
\end{equation*}
It induces the sub-Riemannian metric on the right-invariant distribution. The geodesic flow, which corresponds to the left-invariant sub-Riemannian metric and the horizontal distribution $\mathcal{D}_R$, forms the Hamiltonian system on $T^* H_5$ with the Hamiltonian function:
\begin{align*}
    H(q,\lambda)&=\frac{1}{2\left(x_1^2{+}x_2^2{+}y_1^2{+}y_2^2{+}1\right)}
    \left(\left(x_1^2{+}x_2^2{+}y_2^2{+}1\right) \lambda _1^2{+}\left(x_1^2{+}x_2^2{+}y_1^2{+}1\right) \lambda _2^2{+}\left(x_2^2{+}y_1^2{+}y_2^2{+}1\right) \lambda _3^2 \right.\\
   &{+}\left(x_1^2{+}y_1^2{+}y_2^2{+}1\right) \lambda _4^2{+}\left(x_1^2{+}x_2^2{+}1\right)
   \left(y_1^2{+}y_2^2\right) \lambda _5^2{-}2 x_1 x_2 \lambda _3 \lambda _4{+}2
   \left(y_1^2{+}y_2^2\right) \left(x_1 \lambda _3{+}x_2 \lambda _4\right) \lambda_5\\
   &\left.{+}2 y_1 \left({-}y_2 \lambda _2{+}x_1
   \lambda _3{+}x_2 \lambda _4{+}\left(x_1^2{+}x_2^2{+}1\right) \lambda _5\right) \lambda_1{+}2 y_2 \lambda _2 \left(x_1 \lambda _3{+}x_2 \lambda
   _4{+}\left(x_1^2{+}x_2^2{+}1\right) \lambda _5\right)
    \right).
\end{align*} 

 The Hamiltonian system takes the form:
\begin{align}\label{eq: Hamiltonian system Descartes}
\begin{split}
      \dot{x_1} = \{x_1, H \} &= \frac{(1{+}x_1^2{+}x_2^2{+}y_2^2)\lambda_1{+}y_1(-y_2\lambda_2{+}x_1\lambda_3{+}x_2\lambda_4{+}(1{+}x_1^2{+}x_2^2)\lambda_5)}{1{+}x_1^2{+}x_2^2{+}y_1^2{+}y_2^2},\\
      \dot{x_2} = \{x_2, H \} &= \frac{-y_1y_2\lambda_1(1{+}x_1^2{+}x_2^2)\lambda_2{+}y_1^2\lambda_2{+}y_2(x_1\lambda_3{+}x_2\lambda_4{+}(1{+}x_1^2{+}x_2^2)\lambda_5)}{1{+}x_1^2{+}x_2^2{+}y_1^2{+}y_2^2},\\
      \dot{y_1} = \{y_1, H \} &= \frac{(1{+}x_1^2{+}x_2^2{+}y_2^2)\lambda_3{+}x_1(y_1\lambda_1-x_2\lambda_4{+}y_1^2\lambda_5{+}y_2(\lambda_2{+}y_2\lambda_5))}{1{+}x_1^2{+}x_2^2{+}y_1^2{+}y_2^2},\\
      \dot{y_2} = \{y_2, H \} &= \frac{(1{+}x_1^2{+}x_2^2{+}y_2^2)\lambda_4{+}x_2(y_1\lambda_1-x_1\lambda_3{+}y_1^2\lambda_5{+}y_2(\lambda_2{+}y_2\lambda_5))}{1{+}x_1^2{+}x_2^2{+}y_1^2{+}y_2^2},\\
      \dot{z} = \{z, H \} &= \frac{(1{+}x_1^2{+}x_2^2)y_1\lambda_1{+}(1{+}x_1^2{+}x_2^2)y_2\lambda_2{+}(y_1^2{+}y_2^2)(x_1\lambda_3{+}x_2\lambda_4){+}(1{+}x_1^2{+}x_2^2)(y_1^2{+}y_2^2)\lambda_5}{1{+}x_1^2{+}x_2^2{+}y_1^2{+}y_2^2},\\
      \dot{\lambda_1} = \{\lambda_1, H \} &= - \frac{1}{(1{+}x_1^2{+}x_2^2{+}y_1^2{+}y_2^2)^2}(y_1\lambda_1-x_1\lambda_3-x_2\lambda_4{+}y_1^2\lambda_5{+}y_2(\lambda_2{+}y_2\lambda_5))\cdot\\   
      &\phantom{=-}((1{+}x_2^2{+}y_1^2{+}y_2^2)\lambda_3{+}x_1(y_1\lambda_1-x_2\lambda_4{+}y_1^2\lambda_5{+}y_2(\lambda_2{+}y_2\lambda_5))),\\
      \dot{\lambda_2} = \{\lambda_2, H \} &= - \frac{1}{(1{+}x_1^2{+}x_2^2{+}y_1^2{+}y_2^2)^2}(y_1\lambda_1-x_1\lambda_3-x_2\lambda_4{+}y_1^2\lambda_5{+}y_2(\lambda_2{+}y_2\lambda_5))\cdot\\   
      &\phantom{=-}((1{+}x_1^2{+}y_1^2{+}y_2^2)\lambda_4{+}x_2(y_1\lambda_1-x_1\lambda_3{+}y_1^2\lambda_5{+}y_2(\lambda_2{+}y_2\lambda_5))),\\
      \dot{\lambda_3} = \{\lambda_3, H \} &= \frac{1}{(1{+}x_1^2{+}x_2^2{+}y_1^2{+}y_2^2)^2}(y_1\lambda_1{+}x_4\lambda_2-x_1\lambda_3-x_2\lambda_4-(1{+}x_1^2{+}x_2^2)\lambda_5)\cdot\\
      &\phantom{=-}((1{+}x_1^2{+}x_2^2{+}y_2^2)\lambda_1{+}x_3(-y_2\lambda_2{+}x_1\lambda_3{+}x_2\lambda_4{+}(1{+}x_1^2{+}x_2^2)\lambda_5)),\\
      \dot{\lambda_4} = \{\lambda_4, H \} &= \frac{1}{(1{+}x_1^2{+}x_2^2{+}y_1^2{+}y_2^2)^2}(y_1\lambda_1{+}x_4\lambda_2-x_1\lambda_3-x_2\lambda_4-(1{+}x_1^2{+}x_2^2)\lambda_5)\cdot\\
      &\phantom{=-}(-y_1y_2\lambda_1{+}(1{+}x_1^2{+}x_2^2)\lambda_2{+}y_1^2\lambda_2{+}y_2(x_1\lambda_3{+}x_2\lambda_4{+}(1{+}x_1^2{+}x_2^2)\lambda_5)),\\
      \dot{\lambda_5}=\{ z,H\} &= 0.
    \end{split}
    \end{align}
Here the $z$-coordinate is cyclic, and we can restrict this flow to the level set $\lambda_5=C_0$ and project this restriction onto a 4-dimensional hyperplane $(x_1,x_2,y_1,y_2)$. New system is defined on 8-dimensional symplectic manifold, denoted as $\mathcal{M}_{C_0}$, diffeomorphic to the cotangent bundle of the 4-dimensional hyperplane with the coordinates $(x_1,x_2,y_1,y_2)$. On the manifold $\mathcal{M}_{C_0}$ we introduce the new variables:
 \begin{align*}
 p_1 = \lambda_1 + C_0 y_1, \quad p_2 = \lambda_2 + C_0 y_2, \quad p_3 = \lambda_3,\quad p_4 = \lambda_4,
 \end{align*}
 to express the reduced Hamiltonian function in the form:
\begin{align*}
     H_{C_0}(q,p)=\frac{1}{2} \left(p_1^2 + p_2^2 + p_3^2 + p_4^2 - \frac{(p_3 x_1 + p_4 x_2 - p_1 y_1 - p_2 y_2)^2}{1 + x_1^2 + x_2^2 + y_1^2 + y_2^2}\right).
 \end{align*}
 The Poisson structure on the manifold $\mathcal{M}_{C_0}$ is no longer canonical for $C_0\neq 0$:
\begin{align*}
  \{x_k,p_j\}&=\{y_k,p_{j+2}\}=\delta_{kj},\quad \{p_k,p_{j+2}\}=\delta_{kj}C_0,\\
   \{x_k, y_j\}&=\{x_k,p_{j+2}\}=\{y_k,p_j\}=\{p_j,p_k\}=\{p_{j+2},p_{k+2}\}=0,
\end{align*}
for $j, k=1,2$. Set $W=\frac{(p_3 x_1 + p_4 x_2 - p_1 y_1 - p_2 y_2)}{1 + x_1^2 + x_2^2 + y_1^2 + y_2^2}$ , the corresponding Hamiltonian system is:
\begin{align}\label{eq:system_Ham_W}
\begin{split}
      \dot{x_1} &= \{x_1, H_{C_0} \} = p_1+y_1 W,\\
      \dot{x_2} &= \{x_2, H_{C_0} \} = p_2+y_2 W,\\
      \dot{y_1} &= \{y_1, H_{C_0} \} = p_3-x_1 W,\\
      \dot{y_2} &= \{y_2, H_{C_0} \} = p_4-x_2 W,\\
      \dot{p_1} &= \{p_1, H_{C_0} \} =\left(C_0+W\right) (p_3-x_1 W),\\
      \dot{p_2} &= \{p_2, H_{C_0} \} =\left(C_0+W\right) (p_4-x_2 W),\\
      \dot{p_3} &= \{p_3, H_{C_0} \} = -\left(C_0+W\right) (p_1+y_1 W),\\
      \dot{p_4} &= \{p_4, H_{C_0} \} = -\left(C_0+W\right) (p_2+y_2 W).    
\end{split}
    \end{align}

Change of coordinates to hyper-spherical coordinates $(r,\theta_1, \theta_2, \theta_3)$:
  \begin{align*}
 x_1=r \cos\theta_1 \cos\theta_2, \quad y_1=r \cos\theta_1 \sin\theta_2,\quad
 x_2=r \sin\theta_1 \cos\theta_3, \quad y_2=r \sin\theta_1 \sin\theta_3,
 \end{align*}
with $r\geq 0$, $\theta_1\in\left(0,\frac{\pi}{2}\right)$, $\theta_2, \theta_3\in[0,2\pi)$, and the transformation of the momenta given by:
\begin{align*}
    p_r &=\frac{1}{\sqrt{x_1^2+x_2^2+y_1^2+y_2^2}}( p_1x_1+p_2x_2+p_3y_1+p_4y_2)\\
    p_{\theta_1} &=(p_2x_2+p_4y_2)\sqrt{\frac{x_1^2+y_1^2}{x_2^2+y_2^2}}-(p_1x_1+p_3y_1)\sqrt{\frac{x_2^2+y_2^2}{x_1^2+y_1^2}}\\
    p_{\theta_2} &= p_3x_1-p_1y_1\\
    p_{\theta_3} &= p_4x_2  -p_2y_2, 
\end{align*}

leads to a Hamiltonian function for the reduced system:
\begin{align}\label{eq:Hamiltonian}
 \tilde{H}_{C_0}(r,\theta_1,\theta_2,\theta_3, p_r, p_{\theta_1}, p_{\theta_2}, p_{\theta_3})=\frac{1}{2}\left[p_r^2+\frac{1}{r^2} p_{\theta_1}^2+\frac{1}{r^2}\left(\frac{p_{\theta_2}^2}{\cos^2\theta_1}+\frac{p_{\theta_3}^2}{\sin^2\theta_1}\right)-\frac{1}{1+r^2}(p_{\theta_2}+p_{\theta_3})^2\right].
\end{align}
The Poisson structure on the manifold $\mathcal{M}_{C_0}$ in the hyper-spherical coordinates is:
    \begin{align*}
    \begin{split}
      \{r, p_r\} &=1,\quad \{r,  p_{\theta_k}\}=\{\theta_j, p_r\}= 0,\quad \{\theta_j, p_{\theta_k}\}=\delta_{jk},\quad j,k=1,2,3,\\
      \{p_r,p_{\theta_1}\}&=0, \quad\{p_r,p_{\theta_2}\}=  r C_0\cos^2\theta_1,\quad \{p_r,p_{\theta_3}\}=r C_0 \sin^2\theta_1,\\
       \{p_{\theta_1},p_{\theta_2}\}&=-\{p_{\theta_1},p_{\theta_3}\}=-  r^2 C_0\sin\theta_1\cos\theta_1,\quad \{p_{\theta_2}, p_{\theta_3}\} = 0.
    \end{split}
    \end{align*}
The above structure is non canonical for $C_0\neq 0$ and corresponds to the Hamiltonian system, which describes the motion of a charged particle in a constant magnetic field.

The Hamilton equations are given by
    \begin{align}\label{eq:system_Ham}
    \begin{split}
        \dot{r} &= \{r, \tilde{H}_{C_0}\} = p_r,\\
       \dot{\theta}_1&=\{\theta_1, \tilde{H}_{C_0}\} = \frac{p_{\theta_1}}{r^2},\\
       \dot{\theta}_2&=\{\theta_2, \tilde{H}_{C_0}\} = \frac{p_{\theta_2}}{r^2\cos^2\theta_1}-\frac{p_{\theta_2}+p_{\theta_3}}{1+r^2},\\
       \dot{\theta}_3&=\{\theta_3, \tilde{H}_{C_0}\} = \frac{p_{\theta_3}}{r^2\sin^2\theta_1}-\frac{p_{\theta_2}+p_{\theta_3}}{1+r^2},\\
       \dot{p}_r&=\{p_r, \tilde{H}_{C_0}\}=\frac{1}{r^3}p_{\theta_1}^2+\frac{p_{\theta_2}+p_{\theta_3}}{r(1+r^2)^2}\left((1+r^2)C_0-r^2(p_{\theta_2}+p_{\theta_3})\right)
       +\frac{1}{r^3}\left(\frac{p_{\theta_2}}{\cos^2\theta_1}+\frac{p_{\theta_3}}{\sin^2\theta_1}\right),\\
       \dot{p}_{\theta_1}&=\{p_{\theta_1}, \tilde{H}_{C_0}\}= \frac{1}{r^2}\left(p_{\theta_3}\left(r^2C_0+\frac{p_{\theta_3}}{\sin^2\theta_1}\right)\cot\theta_1 -p_{\theta_2}\left(r^2C_0+\frac{p_{\theta_2}}{\cos^2\theta_1}\right)\tan\theta_1\right),\\
       \dot{p}_{\theta_2}&=\{p_{\theta_2}, \tilde{H}_{C_0}\}=-C_0\cos\theta_1(r p_r\cos\theta_1-p_{\theta_1}\sin\theta_1),\\
       \dot{p}_{\theta_3}&=\{p_{\theta_3}, \tilde{H}_{C_0}\}=-C_0\sin\theta_1(r p_r\sin\theta_1+p_{\theta_1}\cos\theta_1).
    \end{split}
    \end{align}
The flow has four first integrals:
    \begin{align}\label{eq: first_integrals}
    \begin{split}
        \tilde I_1 &= p_{\theta_1}^2+{p_{\theta_2}^2}{\cos^{-2}\theta_1}+{p_{\theta_3}^2}{\sin^{-2}\theta_1}-(p_{\theta_2}+p_{\theta_3})^2,\\
        \tilde I_2 &= p_{\theta_2} + C_0\frac{r^2}{2}  \cos^2\theta_1,\quad \tilde I_3 = p_{\theta_3} + C_0\frac{r^2}{2} \sin^2\theta_1,   \quad \tilde I_4 =\tilde{H}_{C_0}.
    \end{split}
    \end{align}
These integrals are involutive and functionally independent almost everywhere, which implies the integrability of the reduced system in the classical sense. In the cooridinates $(x,y,p)=(x_1,x_2,y_1,y_2,p_1,p_2,p_3,p_4)$ on the manifold $\mathcal{M}_{C_0}$, the first integrals~\eqref{eq: first_integrals} take the form:
\begin{align}\label{eq: first_integrals_descartes}
    \begin{split}
     I_1 =&(y_2p_1-x_1p_4+x_2p_3-y_1p_2)^2+(x_1p_2-x_2p_1-y_1p_4+y_2p_3)^2, \\
     I_2 =& p_3x_1-p_1y_1 + \frac{C_0}{2}(x_1^2+y_1^2), \\
     I_3 =& p_4x_2-p_2y_2 + \frac{C_0}{2}(x_2^2+y_2^2),\\
     I_4 =& H_{C_0}.
     \end{split}
\end{align}

The main objective of this paper is to analyse the
solutions of the reduced Hamiltonian system~\eqref{eq:system_Ham}.

\begin{remark}\label{rmk: origin}
    It should be noted that $r=0$ is merely coordinate singularity of the reduced Hamiltonian system in the hyper-spherical coordinates. Indeed, from the first integrals~\eqref{eq: first_integrals_descartes} in the coordinates $(x,y,p)$  it follows that the solutions of the reduced Hamiltonian system through $(0,0,0,0)$ are characterized by $I_1=I_2=I_3=0$. 
\end{remark}

\begin{remark}\label{rmk: z lift}
Solutions of the reduced Hamiltonian system~\eqref{eq:system_Ham} naturally lift to solutions of the Hamiltonian system~\eqref{eq: Hamiltonian system Descartes}, thus yielding sub-Riemannian geodesics in $\H$. The coordinate $z(t)$ is reconstructed by integrating the corresponding equation in~\eqref{eq: Hamiltonian system Descartes}
\begin{equation}\label{eq: z}
    \dot{z}=\dot{x}_1y_1+\dot{x}_2y_2.
\end{equation}
Geometrical meaning to the evolution of $z(t)$ along sub-Riemannian geodesic line $\gamma: I \rightarrow \H$ starting from $P_0 = \gamma(0) = (x_1(0),x_2(0),y_1(0),y_2(0),z(0))$ can be given together with corresponding formula. First, for $i\in\{1,2\}$, let $\pi_i:\H\rightarrow\mathbb{R}^2$ be projection to subspace $(x_i,y_i)$ and let $(x_i(0),y_i(0))=\pi_i(P_0)$ be projected starting point. Let $A_i(t)$ be the signed sectorial area swept out by position vector from $\pi_i(\gamma(0))$ to $\pi_i(\gamma(t))$, i.e.
\begin{align}\label{eq: sectorial}
    A_i(t)=\frac{1}{2}\int_0^t(x_i(s)\dot{y}_i(s)-\dot{x}_i(s)y_i(s))ds.
\end{align}
From identity 
\begin{align*}
    \frac{d}{dt}(x_iy_i)=\dot{x}_iy_i+x_i\dot{y}_i
\end{align*}
we have 
\begin{align}\label{eq: total derivative}
    x_i(t)y_i(t)-x_i(0)y_i(0)=\int_0^t(\dot{x}_i(s)y_i(s)+x_i(s)\dot{y}_i(s))ds.
\end{align}

From~\eqref{eq: sectorial} and~\eqref{eq: total derivative} it follows that 
\begin{align*}
    \int_0^t\dot{x}_i(s)y_i(s)s=\frac{1}{2}\big(x_i(t)y_i(t)-x_i(0)y_i(0))-A_i(t).
\end{align*}
This yields formulae for $z(t)$ of the form:
\begin{align*}
    z(t)&=z(0)+\frac{1}{2}(x_1(t)y_1(t)+x_2(t)y_2(t)-x_1(0)y_1(0)-x_2(0)y_2(0)\big)- A_1(t)-A_2(t).
\end{align*}

We illustrate the behaviour of the solutions through several numerical examples presented at the end of the paper.

\end{remark}

\section{Classification of geodesics}\label{sec3}

We are going to describe solutions that belong to the 4-dimensional invariant manifold. Since $ \tilde I_1, \tilde I_2, \tilde I_3$ and $\tilde I_4$, defined in~\eqref{eq: first_integrals}, are first integrals, there are constants $C_1, C_2, C_3, C_4 \in \mathbb{R}$ such that 
\begin{align*}
    \tilde I_1=C_1,\quad \tilde I_2=C_2,\quad \tilde I_3 = C_3, \quad \tilde I_4=C_4.
\end{align*} 
The values of the first integrals determine the initial points from which the corresponding geodesics start; in particular, for the geodesic through the origin all values $C_1, C_2, C_3$ vanish. In the analysis that follows, we study the first-integral constants in their generic form, thereby encompassing all possible choices of their values, including the special case in which they all vanish.

With~\eqref{eq: first_integrals} we derive
\begin{gather}\label{eq:p_teta_12}
    \begin{split}
    &p_{\theta_2}+p_{\theta_3} = C_2+C_3-C_0\frac{r^2}{2},\\
    &\frac{C_1}{r^2}=\frac{1}{r^2}(p_{\theta_1}^2+{p_{\theta_2}^2}{\cos^{-2}\theta_1}+{p_{\theta_3}^2}{\sin^{-2}\theta_1}-(p_{\theta_2}+p_{\theta_3})^2), 
    \end{split}  
\end{gather}
By further using the conservation laws and the expressions~\eqref{eq:Hamiltonian},~\eqref{eq: first_integrals} and~\eqref{eq:p_teta_12} we can reduce the problem to a problem with one degree of freedom for the variable $r$:
\begin{align}\label{eq: Doted r and C4}
     2C_4=\Dot{r}^2+\frac{1}{r^2} C_1+\frac{1}{r^2(1+r^2)}(C_2+C_3-C_0\frac{r^2}{2})^2.
\end{align}
For $C_4\neq 0$, equation~\eqref{eq: Doted r and C4} in a more compact form is:
\begin{align}
 \Dot{r}^2&=2 C_4 \frac{f(r)}{r^2(1+r^2)},
\end{align}
where $f(r)$ is a polynomial of degree 4, in the general case:
\begin{align}\label{eq: f(r)}
    \begin{split}
    f(r)&=A r^4+ B r^2 + C,\\
    A = 1-\frac{C_0^2}{8C_4}, \quad B =1-\frac{C_1}{2C_4}&+\frac{(C_2+C_3)C_0}{2C_4}, \quad 
    C = -\frac{C_1+(C_2+C_3)^2}{2C_4}.    
    \end{split}  
\end{align}
We can restrict the discussion to the case where the solutions are parametrised by the arc length parameter, which is equivalent to $C_4=\frac{1}{2}$ according to the Cauchy-Schwarz integral inequality between length and energy function, equivalent to $C_4=\frac{1}{2}$. This simplifies the problem to:
\begin{align}\label{eq:Doted_r}
    \begin{split}
    \Dot{r}^2&= \frac{f(r)}{r^2(1+r^2)},\quad \textrm{i.e.} \quad\Dot{r}= \pm\frac{\sqrt{f(r)}}{r\sqrt{1+r^2}}, \\
    A = 1-\frac{C_0^2}{4}, \quad B &=1-C_1+(C_2+C_3)C_0, \quad 
    C = -C_1-(C_2+C_3)^2.  
    \end{split}  
\end{align}
With~\eqref{eq: first_integrals} and the Hamilton's equations we have:
\begin{align}
    C_1=r^4\dot{\theta}_1^2+\frac{C_2^2}{{\cos{\theta_1}}^2}+\frac{C_3^2}{{\sin{\theta_1}}^2}-(C_2+C_3)^2,\notag\\
    \dot{\theta}_1=\pm\frac{1}{r^2}\sqrt{C_1+(C_2+C_3)^2-\frac{C_2^2}{{\cos^2{\theta_1}}}-\frac{C_3^2}{{\sin^2{\theta_1}}}}.\label{eq:teta_prim_2}
\end{align}
The conservation laws leads to the following equations for the remaining coordinates:
\begin{align}\label{eq:theta_23}
\begin{split}
    \dot \theta_2 = \frac{C_2}{r^2\cos^2{\theta_1}}-\frac{C_2+C_3}{1+r^2}-\frac{C_0}{2(1+r^2)},\\
   \dot \theta_3 = \frac{C_3}{r^2\sin^2{\theta_1}}-\frac{C_2+C_3}{1+r^2}-\frac{C_0}{2(1+r^2)}.
\end{split}
\end{align}

\begin{proposition}\label{thm:quadratures}
The magnetic geodesic flow~\eqref{eq:system_Ham} is solvable by
quadratures.
\end{proposition}
\begin{proof}
    The integration of the entire problem can be carried out in four steps. First, integrate inverted~\eqref{eq:Doted_r} to implicitly obtain $r=r(t)$:
\begin{align}\label{eq:Doted_t}
    \Dot{t}= \pm\frac{r\sqrt{r^2+1}}{\sqrt{f(r)}},  
\end{align}
\begin{equation*}
 t+t_0=\pm\int\frac{r\sqrt{r^2+1}}{\sqrt{f(r)}}dr, \quad t_0\in\mathbb{R}.
\end{equation*}
After the initial choice of sign, the expression on the right-hand side of~\eqref{eq:Doted_t} has a constant sign, which means that $t(r)$ is strictly monotonic and therefore allows inversion.

The chain rule $\frac{d\theta_1}{dt}=\frac{d\theta_1}{dr}\frac{dr}{dt}$ allows the integration~\eqref{eq:teta_prim_2} to implicitly define $\theta_1(r)$. It should be noted that we can obtain $\theta_1(r)$ without knowing $r(t)$, together with the previous step $\theta_1(t)=\theta_1(r(t))$:

\begin{align*}
 \int\frac{d\theta_1}{\sqrt{C_1+(C_2+C_3)^2-\frac{C_2^2}{{\cos^2{\theta_1}}}-\frac{C_3^2}{{\sin^2{\theta_1}}}}}&=\pm\int\frac{\sqrt{1+r^2}}{r\sqrt{f(r)}}dr.
\end{align*}
Using the chain rules $\frac{d\theta_2}{dt}=\frac{d\theta_2}{dr}\frac{dr}{dt}$ and $\frac{d\theta_3}{dt}=\frac{d\theta_3}{dr}\frac{dr}{dt}$, the equations~\eqref{eq:theta_23} and the results of the previous steps, the remaining coordinate functions are recovered by :
\begin{align*}
 \theta_2+\theta_2^0 &= \int \frac{r\sqrt{r^2+1}}{\sqrt{f(r)}}\left(\frac{C_2}{r^2\cos^2{\theta_1(r)}}-\frac{C_2+C_3+\frac{1}{2}C_0}{1+r^2}\right)dr, \quad \theta_2^0\in \mathbb{R};\\
 \theta_3 + \theta_3^0 &= \int \frac{r\sqrt{r^2+1}}{\sqrt{f(r)}}\left(\frac{C_3}{r^2\sin^2{\theta_1(r)}}-\frac{C_2+C_3+\frac{1}{2}C_0}{1+r^2}\right)dr, \quad \theta_3^0\in \mathbb{R}.
\end{align*}
Note that the second step is necessary for the last two and that we can obtain $\theta_i=\theta_i(r),\ i\in\{2,3\}$ without integrating the first step, while to define $\theta_1=\theta_1(t)$ the integration of~\eqref{eq:Doted_r} is necessary.
\end{proof}
The following statements will be useful later.
\begin{lemma}\label{lemma:C1}
     $C_1 = p_{\theta_1}^2+{p_{\theta_2}^2}{\cos^{-2}\theta_1}+{p_{\theta_3}^2}{\sin^{-2}\theta_1}-(p_{\theta_2}+p_{\theta_3})^2\geq 0$.
\end{lemma}
\begin{proof}
    Statement of lemma follows from the inequality
 \begin{equation}\label{eq:aux_C1_ineq}
 \frac{x^2}{\cos^2{\theta}}+\frac{y^2}{\sin^2{\theta}}\geq (x+y)^2
 \end{equation}
 applied to $x=p_{\theta_2}, y=p_{\theta_3}$ and $\theta=\theta_1$. To prove  inequality~\eqref{eq:aux_C1_ineq} we first notice that for $x$ and $y$ with different signs, inequality is trivial. For $\sgn(x)=\sgn(y)$, we define $f:(0,\frac{\pi}{2})\rightarrow \mathbb{R}, f(\theta)=\frac{x^2}{\cos^2{\theta}}+\frac{y^2}{\sin^2{\theta}}$ and find the minima for $f$. Since $f'(\theta)=2\sin{\theta}\cos{\theta}(\frac{x^2}{\cos^4{\theta}}-\frac{y^2}{\sin^4{\theta}})$, the minimum at $\theta_{min}$ is such that $\tan^4{\theta_{min}}=\frac{y^2}{x^2}$. Now~\eqref{eq:aux_C1_ineq} is equivalent to
 \begin{gather*}
 f(\theta_{min})\geq(x+y)^2 \Leftrightarrow x^2\tan^2{(\theta_{min})}+y^2\ctan^2{(\theta_{min})}\geq 2xy \Leftrightarrow x^2\frac{y}{x}+y^2\frac{x}{y}\geq 2xy.
 \end{gather*}
\end{proof}
\begin{cor}\label{cor:C_lessthan_0}
  The parameter $C$ is non-positive, that is, $C=-C_1-(C_2+C_3)^2\leq0$. Furthermore, $C=0$ implies $B=1$.
\end{cor}
\begin{remark}
We note here that it is not possible for $f(r)$ to be a zero polynomial, which means that for $C_4\neq 0$ $r(t)$ is never a constant function. Moreover, $A\leq 1$.
\end{remark}

Let us notice that if $C_4=0$ in~\eqref{eq: Doted r and C4}, motion degenerates to a point. This assertion follows from the following: Using Lemma~\ref{lemma:C1} with $0=C_4=\tilde{H}_C$ implies that all the moments $p_r, p_{\theta_1}, p_{\theta_2}$ and $p_{\theta_3}$ are equal to zero. Now, only solution to~\eqref{eq:system_Ham} are constants. In the sequel we assume $C_4\neq 0$.

We fix values of first integrals by: $\tilde I_1=C_1$, $\tilde I_2=C_2$, $\tilde I_3 = C_3$, $\tilde I_4=C_4$ and formulate our main result as follows.
\begin{theorem}\label{thm: classification reduced}
 The reduced Hamiltonian system~\eqref{eq:system_Ham} has trajectories of two types:
 \begin{enumerate}[label=(\roman*)]

\item\label{type1} Trajectories of the first type satisfy $r(t)\ge r_0$ for all $t$, where $r_0$ is the largest non-negative root of~\eqref{eq: f(r)}.

\item\label{type2} Trajectories of the second type satisfy $r_1\le r(t)\le r_2$ for all $t$, where $r_1<r_2$ are the positive roots of~\eqref{eq: f(r)}.
 \end{enumerate}
\end{theorem}
\begin{remark}
From~\eqref{eq:Doted_r} it is seen that motion is possible only in regions where $f(r)>0$ (see Figure~\ref{fig:fig1} and Figure~\ref{fig:fig2}). This implies a strong dependence between the constants of motion $C_k$, $k\in\{1,2,3,4\}$, and the set of admissible points. More precisely, given constants $C_k$, the motion is confined to certain regions and, consequently, not all initial conditions are admissible; conversely, the choice of initial conditions restricts the possible values of the constants of motion. For example, for geodesics passing through the origin holds $C_1=C_2=C_3=0$.

Theorem~\ref{thm: classification reduced} provides a lower bound for the radial distance along a trajectory; however, this bound depends on the constants of motion and, in general, may be arbitrarily small.
\end{remark}

From the Remark~\ref{rmk: z lift} it follows that the classification in Theorem~\ref{thm: classification reduced} extends naturally to the classification of the original Hamiltonian system~\eqref{eq: Hamiltonian system Descartes} in $\H$. 
\begin{cor}
    Every geodesic of the reduced system~\eqref{eq:system_Ham_W} lifts to a geodesic of the system~\eqref{eq: Hamiltonian system Descartes} via the reconstruction equation~\eqref{eq: z}.
\end{cor}

In the following, we prove Theorem~\ref{thm: classification reduced}.

\subsection{Case: $A=0$}
\begin{figure}[H]
\centering
      {\includegraphics[width=.3\textwidth]{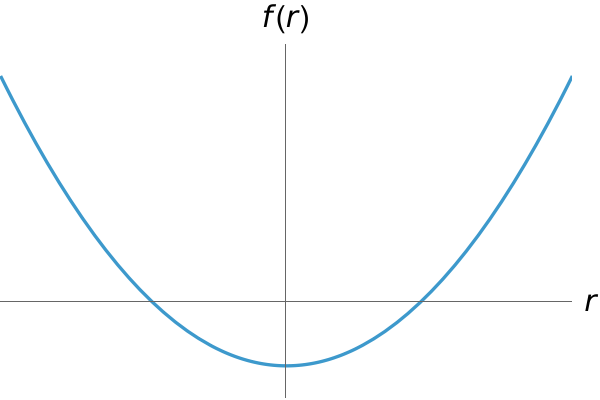}}
  \caption{$B>0, C\leq0$.}
  \label{fig:fig1}
\end{figure}

In this case $C_0=2$ and $f(r)$ simplify to:
\begin{align*}
    f(r)&=B r^2+C.
\end{align*}
Since $f(r)\geq 0$, only $B>0$ is possible in this case (shown in Figure~\ref{fig:fig1}) and motion is only possible in the region $r > r_0=\sqrt{-\frac{C}{B}}$.
Now we continue with~\eqref{eq: f(r)} to implicitly find $r(t)$:
\begin{gather*}
    \Dot{r}=\pm \frac{\sqrt{Br^2+C}}{r\sqrt{1+r^2}}\\
    t = \pm \frac{1}{\sqrt{B}}\int \frac{r\sqrt{1+r^2} dr}{\sqrt{r^2+C/B}}.
\end{gather*} 
With $a=-\frac{C}{B}\geq0$, the substitutions $1+r^2=z^2$ and $p=\frac{z}{\sqrt{a+1}}$ we solve the integral on the right-hand side as follows:
\begin{align*}
    \pm \frac{1}{\sqrt{B}}\int \frac{\sqrt{1+r^2}r dr}{\sqrt{r^2+C/B}}&=\pm \frac{1}{\sqrt{B}}\int \frac{z^2 dz}{\sqrt{z^2-(a+1)}}\\
    &=\pm \frac{1}{\sqrt{B}}\left(\int \sqrt{z^2-(a+1)}dz + (a+1)\int \frac{dz}{\sqrt{z^2-(a+1)}}\right)\\
    &=\pm \frac{a+1}{\sqrt{B}}\left(\int \sqrt{p^2-1}dp + \int \frac{dp}{\sqrt{p^2-1}}\right)\\
    &=\pm \frac{a+1}{\sqrt{B}}\left(\frac{1}{2}p\sqrt{p^2-1}-\frac{1}{2}\ln(p+\sqrt{p^2-1})+\ln(p+\sqrt{p^2-1})\right)\\
    &=\pm \frac{1+a}{\sqrt{B}}\left(\frac{1}{2}\sqrt{\frac{1+r^2}{1+a}}\sqrt{\frac{r^2-a}{1+a}}+\frac{1}{2}\ln(\sqrt{\frac{1+r^2}{1+a}}+\sqrt{\frac{r^2-a}{1+a}})\right).
\end{align*} 
Finally:
\begin{align*}
    t=t_0\pm \frac{1}{2\sqrt{B}}\left(\sqrt{(1+r^2)(r^2-a)}+(1+a)\ln\big(\sqrt{1+r^2}+\sqrt{r^2-a}\big)-\frac{1+a}{2}\ln(1+a)\right).
\end{align*} 

\begin{example} Let us consider example where $A=0, B=-C$ and plot $r(t)$ and $\theta_1(t)$.
\begin{figure}[H]
\centering
      {\includegraphics[width=.3\textwidth]{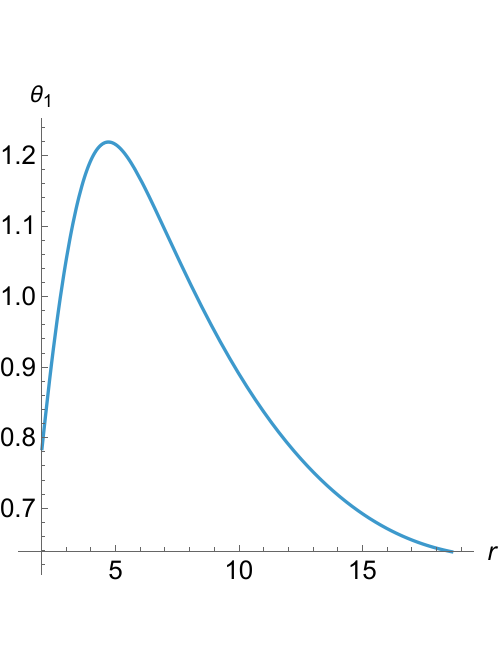}}
      \quad
      {\includegraphics[width=.3\textwidth]{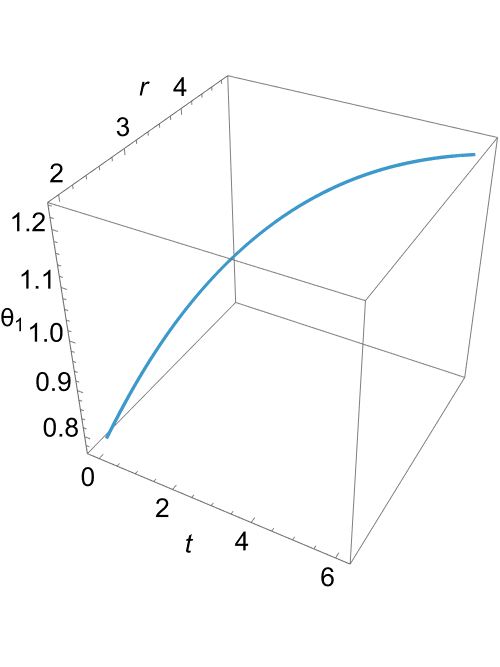}}
      \quad
      {\includegraphics[width=.3\textwidth]{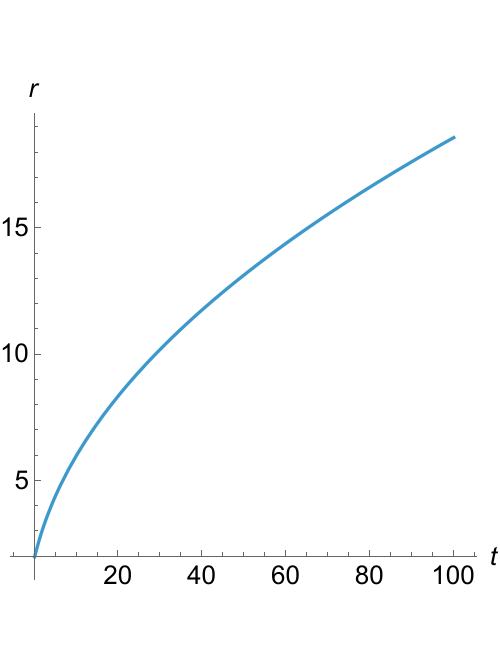}}
  \caption{Plots of $(r(t),\theta_1(t))$ (left), $(t,r(t),\theta_1(t))$ (centre) and $(t,r(t))$ (right) for $f(r)=\frac{9}{2}r^2-\frac{9}{2}$.}
  \label{fig:example}
\end{figure}
\end{example}

\subsection{Case: $A\neq0$}
\begin{figure}[H]
\centering
      {\includegraphics[width=.3\textwidth]{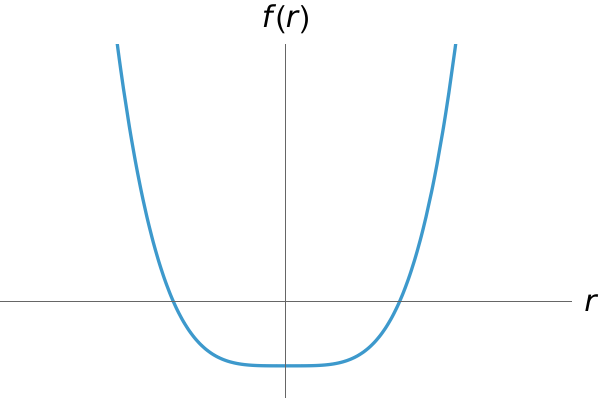}}
      \quad
      {\includegraphics[width=.3\textwidth]{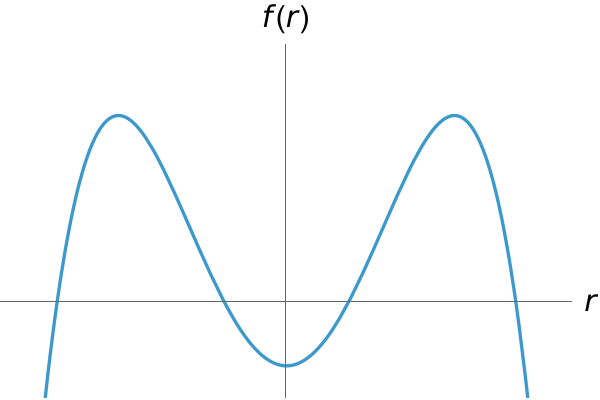}}
  \caption {Case $A>0$ (on the left) and Case $A<0$ (on the right).}
  \label{fig:fig2}
\end{figure}

For the case $A\neq0$ we have two sub-cases in which many situations are eliminated thanks to Lemma~\ref{cor:C_lessthan_0}. The cases of interest, which will be explained later, are shown in Figure~\ref{fig:fig2}.

\subsubsection{Case: $A>0$}

If we look at the discriminant $D=B^2-4AC$ of the quadratic equation $Ax^2+Bx+C$, we see that it is always non-negative, since $A>0$ and $C\leq0$. This means that we can write $f(r)=A(r^2-\bar\alpha)(r^2-\bar\beta)$, with $C=\frac{\bar\alpha\bar\beta}{A}\leq0$, and that we have exactly two real roots $\bar\alpha\geq0$ and $\bar\beta\leq0$.
The function $f(r)=A(r+\alpha)(r-\alpha)(r^2+\beta^2)$ is non-negative for $r\geq\alpha$.
We distinguish between two cases: $C=0$ and $C<0$.

In the first case, we have that the equation~\eqref{eq:Doted_r} reduces to:
\begin{align}\label{eq: t doted}
    \dot{t}^2=\frac{1+r^2}{Ar^2+1}.
\end{align}
Note that if $C_0=0$, i.e. $A=1$ we have simple form for $t(r)$, namely $t(r)=t_0\pm r$.

In general case, $C_0\neq 0$, we integrate~\eqref{eq: t doted} as
\begin{align*}
    t=t_0\pm\int \frac{1+r^2}{\sqrt{(1+r^2)(Ar^2+1)}}dr,
\end{align*}
which solves to
\begin{align*}
    t=t_0\pm F(\phi,k)\pm\frac{1}{A}\left(\dn{u}\sc1{u}-E(u)\right), \quad \phi=\am{u}, \sc1{u}=r, k=\sqrt{1-A^2}.
\end{align*}
In this case $\theta_1$ is constant function, and $\theta_2$ and $\theta_3$ are linear functions as stated in the Corollary~\ref{cor:c0}.

In the second case, if $C<0$, we solve the equation:
\begin{align}\label{eq:tr}
\dot{t}^2=\frac{r^2(1+r^2)}{A(r+\alpha)(r-\alpha)(r^2+\beta^2)}.
\end{align}
We must first separately consider the case where $\beta=1$. The solution of~\eqref{eq:tr} is:
\begin{align*}
    t(r)=t_0\pm\sqrt{\frac{r^2-\alpha^2}{A}}.
\end{align*}

If $\beta\neq 1$ substitution $1+r^2=z^2$ reduces the problem to an elliptic integral of the second kind.
\begin{align*}
 t&=\pm\frac{1}{\sqrt{A}}\int\frac{r\sqrt{1+r^2}dr}{\sqrt{(r+\alpha)(r-\alpha)(r^2+\beta^2)}}=\pm\frac{1}{\sqrt{A}}\int\frac{z^2dz}{\sqrt{(z^2-(\alpha^2+1))(\beta^2-1+z^2)}}.
\end{align*}
We distinguish between cases where $\beta^2>1$ and $\beta^2<1$.

In the case of $\beta^2<1$ integral is further solved with substitution $\frac{z}{\sqrt{\alpha^2+1}}=\ns{u}$ (see \cite[p. 192]{BF1971}) as:
\begin{align*}    
    t&=\pm\frac{1}{\sqrt{A}}\int\frac{z^2dz}{\sqrt{(z^2-(\alpha^2+1))(\beta^2-1+z^2)}}=\pm\frac{1}{\sqrt{A}}\int\frac{\sqrt{\alpha^2+1}(\frac{z}{\sqrt{\alpha^2+1}})^2d(\frac{z}{\sqrt{\alpha^2+1}})}{\sqrt{((\frac{z}{\sqrt{\alpha^2+1}})^2-1)(\frac{\beta^2-1}{\alpha^2+1}+(\frac{z}{\sqrt{\alpha^2+1}})^2)}}\\
    &=\pm\frac{\sqrt{\alpha^2+1}}{\sqrt{A}}\int \ns^2{u}du+t_0=\pm\frac{\sqrt{\alpha^2+1}}{\sqrt{A}}(u-E(u)-\dn{u}\cs{u})+t_0,
\end{align*}
where $\ns u=\sqrt{\frac{r^2+1}{\alpha^2+1}}, t_0\in\mathbb{R}$.

In case of $\beta^2>1$, the integration is performed with the substitution $\frac{z}{\sqrt{\alpha^2+1}}=\nc{u}$ (see \cite[p. 193]{BF1971}) as:
\begin{align*}
 t&=\pm\frac{1}{\sqrt{A}}\int\frac{z^2dz}{\sqrt{(z^2-(\alpha^2+1))(\beta^2-1+z^2)}}=\pm\frac{1}{\sqrt{A(\alpha^2+\beta^2)}}\int\frac{(\alpha^2+1)(\frac{z}{\sqrt{\alpha^2+1}})^2d(\frac{z}{\sqrt{\alpha^2+1}})}{\sqrt{((\frac{z}{\sqrt{\alpha^2+1}})^2-1)(\frac{\beta^2-1}{\alpha^2+\beta^2}+\frac{\alpha^2+1}{\alpha^2+\beta^2}(\frac{z}{\sqrt{\alpha^2+1}})^2)}}\\
 &=\pm\frac{\alpha^2+1}{\sqrt{A(\alpha^2+\beta^2)}}\int \nc^2{u}du+t_0=\pm\frac{\alpha^2+1}{k'^2\sqrt{A(\alpha^2+\beta^2)}}(k'^2u-E(u)+\dn{u}\sc1{u})+t_0,
\end{align*}
where $\nc u=\sqrt{\frac{r^2+1}{\alpha^2+1}}, k'^2=\frac{\alpha^2+1} {\alpha^2+\beta^2}, t_0\in\mathbb{R}$.

Finally, it is worth investigating the asymptotic behaviour of the trajectories of the first type when $r\rightarrow+\infty$. From~\eqref{eq:teta_prim_2} and~\eqref{eq:theta_23} it can be seen that we have a very simple form of the functions $\theta_i$, $i\in\{1,2,3\}$.
\begin{cor}\label{cor:rinf}
    If $A\geq 0$, when $r\rightarrow+\infty$, $\theta_1$ tends to a constant function and $\theta_2$ and $\theta_3$ tend to linear functions in~$t$.
\end{cor}

\subsubsection{Case: $A<0$}
Here, too, we have two different cases.

In the first, similar to the previous consideration, we find that the function $f(r)$ has four real roots: $\alpha$, $-\alpha$, $\beta$ and $-\beta$. Without loss of generality, we can assume that $\alpha<\beta$. The function $f(r)$ has the form $f(r)=A(r-\alpha)(r+\alpha)(r-\beta)(r+\beta)$ and is non-negative in the interval $[\alpha,\beta]$.
We solve the equation:
\begin{align*}
 \dot{t}^2&=\frac{r^2(1+r^2)}{A(r+\alpha)(r-\alpha)(r+\beta)(r-\beta)},\quad r\in(\alpha,\beta).
\end{align*}
The substitution $1+r^2=z^2$ reduces the problem to an elliptic integral of the second kind, and with the subsequent substitution $z=\sqrt{\alpha^2+1}\nd{u}$ the solution is obtained by the elliptic function $E$:

\begin{align*}
 t&=\pm\frac{1}{\sqrt{|A|}}\int\frac{r\sqrt{1+r^2}dr}{\sqrt{(r+\alpha)(r-\alpha)(\beta+r)(\beta-r)}}=\pm\frac{1}{\sqrt{|A|}}\int\frac{z^2dz}{\sqrt{(z^2-(\alpha^2+1))(\beta^2+1-z^2)}}\\
 &=\pm\frac{\alpha^2+1}{\sqrt{|A|(\beta^2+1)}}\int\frac{(\frac{z}{\sqrt{\alpha^2+1}})^2d(\frac{z}{\sqrt{\alpha^2+1}})}{\sqrt{((\frac{z}{\sqrt{\alpha^2+1}})^2-1)(1-\frac{\alpha^2+1}{\beta^2+1}(\frac{z}{\sqrt{\alpha^2+1}})^2)}}=\pm\frac{\alpha^2+1}{\sqrt{|A|(\beta^2+1)}}\int\nd^2{u}du\\
 &=\pm\frac{\alpha^2+1}{k'^2\sqrt{|A|(\beta^2+1)}}(E(u)-k^2\sn{u}\cd{u})+c_0,
\end{align*}
where $\nd u=\sqrt{\frac{r^2+1}{\alpha^2+1}}, k^2=\frac{\alpha^2+1}{\beta^2+1}, c_0\in\mathbb{R}$.

The second case is that there are no real roots for the function $f(r)$. In this case, however, the function is negative in its domain, so we have to discard it.

This concludes the consideration and proves the Theorem~\ref{thm: classification reduced}. It also proves the following theorem.

\begin{theorem}
The first step for quadrature integration, i.e. the formulas for $r=r(t)$, are defined locally (implicitly or by inverting elliptic integrals) as follows:
\begin{enumerate}[label=\alph*)]
    \item if $A=0$: 
    \begin{align*}
        t=t_0\pm \frac{1+a}{\sqrt{B}}\left(\frac{1}{2}\sqrt{\frac{1+r^2}{1+a}}\sqrt{\frac{r^2-a}{1+a}}+\frac{1}{2}\ln(\sqrt{\frac{1+r^2}{1+a}}+\sqrt{\frac{r^2-a}{1+a}})\right),\quad a=-\frac{C}{B}.
    \end{align*}
    \item if $0<A\leq1, C=0$:
    \begin{align*}
        t=t_0\pm F(\phi,k)\pm\frac{1}{A}\left(\dn{u}\sc1{u}-E(u)\right), \quad \phi=\am{u}, \sc1{u}=r, k=\sqrt{1-A^2}.
    \end{align*}
    \item if $0<A\leq1, C<0$ and $\pm i$ are roots of $f(r)$:  
    \begin{align*}
        t(r)=t_0\pm\sqrt{\frac{r^2-\alpha^2}{A}}.
    \end{align*}
    \item if $0<A\leq1, C<0$ and $\beta^2<1$, where $\beta$ is complex, and $\pm\alpha$ are real roots of $f(r)$:
    \begin{align*}    
    t=\pm\frac{\sqrt{\alpha^2+1}}{\sqrt{A}}(u-E(u)-\dn{u}\cs{u})+t_0,
\end{align*}
    \item if $0<A\leq1, C<0$ and $\beta^2>1$, where $\beta$ is complex root, and $\pm\alpha$ are real roots of $f(r)$:
    \begin{align*}    
    t=\pm\frac{\alpha^2+1}{k'^2\sqrt{A(\alpha^2+\beta^2)}}(k'^2u-E(u)+\dn{u}\sc1{u})+t_0,\quad \nc u=\sqrt{\frac{r^2+1}{\alpha^2+1}}, k'^2=\frac{\alpha^2+1} {\alpha^2+\beta^2}.
    \end{align*}
    \item if $A<0$ and $\pm\alpha, \pm\beta$ are real roots of $f(r)$:
    \begin{align*}    
    t=\pm\frac{\alpha^2+1}{k'^2\sqrt{|A|(\beta^2+1)}}(E(u)-k^2\sn{u}\cd{u})+t_0,\quad \nd u=\sqrt{\frac{r^2+1}{\alpha^2+1}}, k^2=\frac{\alpha^2+1}{\beta^2+1}.
    \end{align*}
\end{enumerate}
\end{theorem}

Let us consider the case where $C=0$, which contains geodesics passing through coordinate origin. We have a similar situation as in Corollary~\ref{cor:rinf}.
\begin{cor}\label{cor:c0}
    If $C=0$, $\theta_1$ is a constant function and $\theta_2$ and $\theta_3$ are linear functions in $t$.
\end{cor}

We conclude the paper with numerical examples of the sub-Riemannian geodesic curves of the original system on $\H$. 
\begin{example} 
In this example numerically obtained sub-Riemannian geodesic curves passing through coordinate origin are  visualized for both types of trajectories from Theorem~\ref{thm: classification reduced} (see Figure~\ref{fig: z(t) A<0} and Figure~\ref{fig: z(t) A>0}). Corresponding function $f(r)$ from~\eqref{eq: f(r)} is of the form $f(r)=-r^4+r^2$ for geodesic in the Figure~\ref{fig: z(t) A<0} and of the form $f(r)=\frac{1}{2}r^4+r^2$ for geodesic in the Figure~\ref{fig: z(t) A>0}. The left and right panels of figures show the dependence of the $z$-coordinate and the radius $r=\sqrt{x_1^2+x_2^2+y_1^2+y_2^2}$ on $t$. The centre panels show the projection of the sub-Riemannian geodesic onto the $(x_1,y_1,z)$ subspace.  
    \begin{figure}[H]
\centering
      {\includegraphics[width=.25\textwidth]{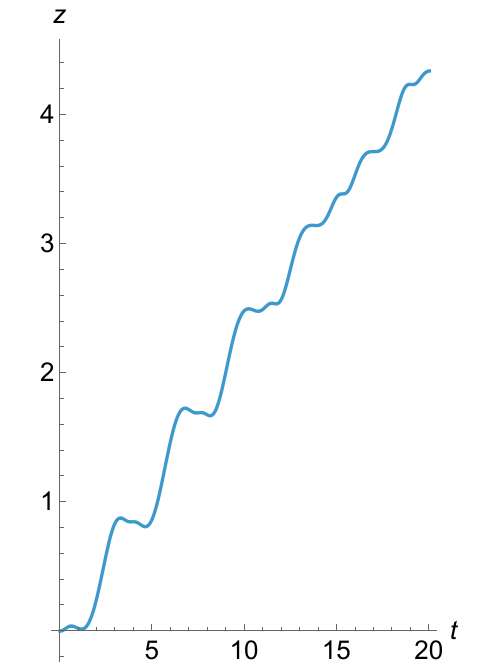}}
      \quad
      {\includegraphics[width=.25\textwidth]{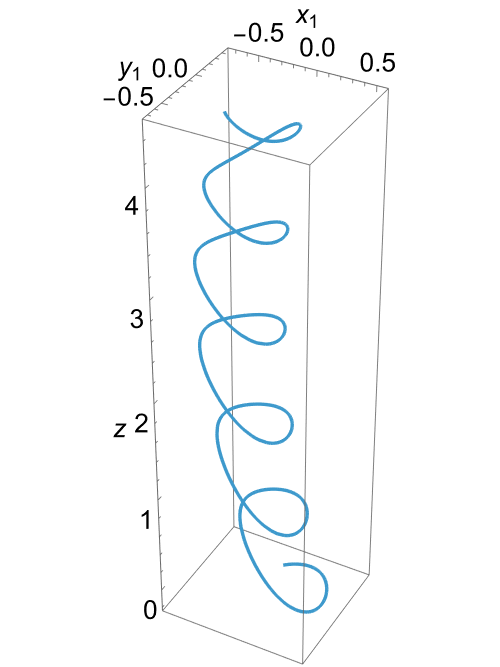}}
      \quad
      {\includegraphics[width=.25\textwidth]{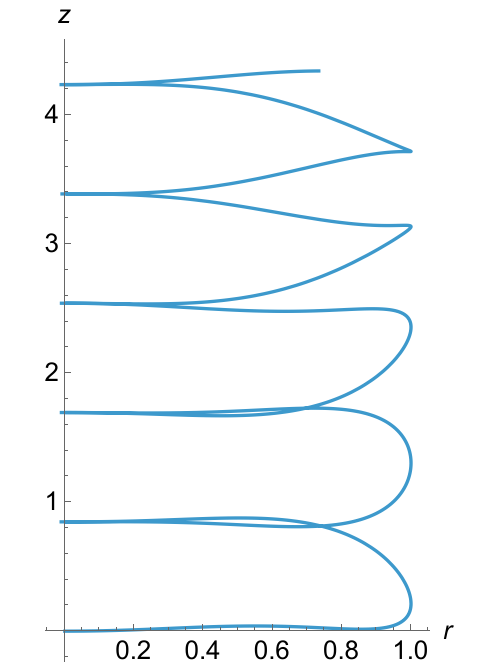}}
  \caption {Plots of $(t,z(t))$ (left), $(r(t),z(t))$ (right) and projection onto the $(x_1,y_1,z)$ subspace (centre) of a sub-Riemannian geodesic of type~\ref{type2} passing through origin.}
  \label{fig: z(t) A<0}
\end{figure}
\end{example}

    \begin{figure}[H]
\centering
      {\includegraphics[width=.25\textwidth]{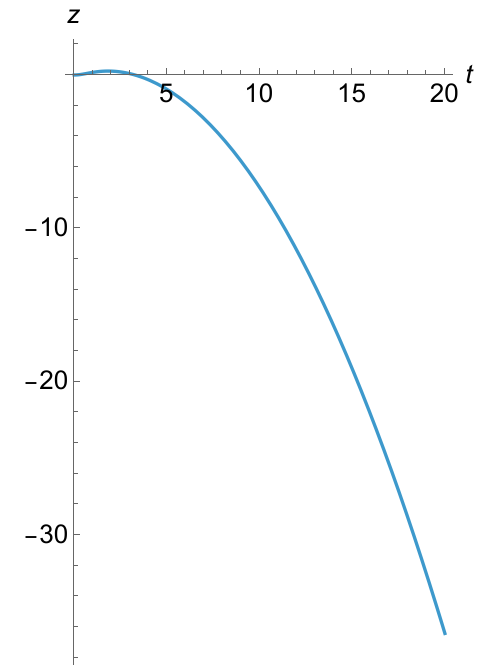}}
      \quad
      {\includegraphics[width=.25\textwidth]{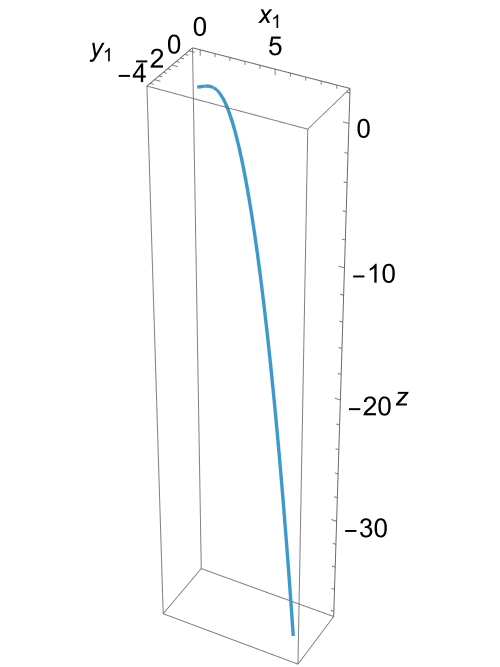}}
      \quad
      {\includegraphics[width=.25\textwidth]{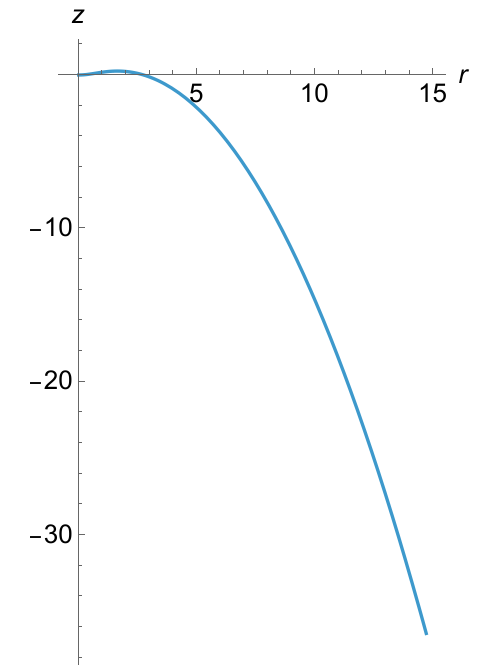}}
  \caption {Plots of $(t,z(t))$ (left), $(r(t),z(t))$ (right) and projection onto the $(x_1,y_1,z)$ subspace (centre) of sub-Riemannian geodesic of type~\ref{type1}  passing through origin.}
  \label{fig: z(t) A>0}
\end{figure}

In the previous example with a geodesic line through origin it was not apparent that projection of a geodesic line on the $(x_1,x_2,y_1,y_2)$ subspace is bounded both bellow and above. Following example demonstrates that property. 
\begin{example} Here we visualize geodesic line through point $(0,-\sqrt{2},-\sqrt{2},0,0)$ on level set $C_0= 2\sqrt{2}$ and with values of first integrals as: $C_1=1, C_2=\sqrt{2}, C_3=\sqrt{2}$ and $C_4=\frac{1}{2}$, resulting in $f(r)=-r^4+8r^2-9$. This example illustrates trajectories of type~\ref{type2} from Theorem~\ref{thm: classification reduced} and they satisfy $4-\sqrt{7}\leq r^2(t)\leq 4+\sqrt{7}$ (see Figure~\ref{fig: ring geodesic}). To visualize objects in five-dimensional real space, we interpret the coordinate pairs $(x_1, y_1)$ and $(x_2, y_2)$ as complex numbers $z_1$ and $z_2$, respectively, and represent their moduli along the first two coordinates of the plot, while retaining the fifth coordinate as the third axis. The resulting three-dimensional plot therefore has coordinates $(\|z_1\|,\|z_2\|, z)$.
    \begin{figure}[H]
\centering
      {\includegraphics[width=.25\textwidth]{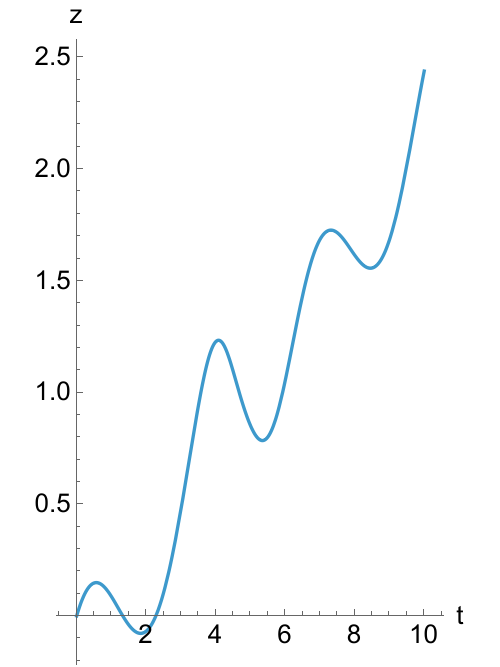}}
      \quad
      {\includegraphics[width=.25\textwidth]{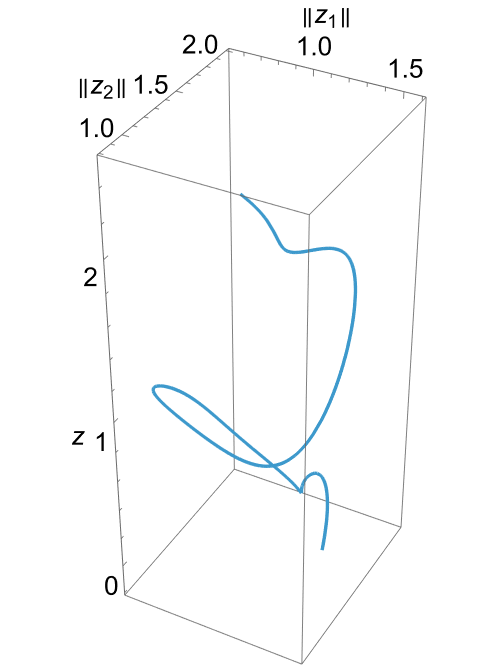}}
      \quad
      {\includegraphics[width=.25\textwidth]{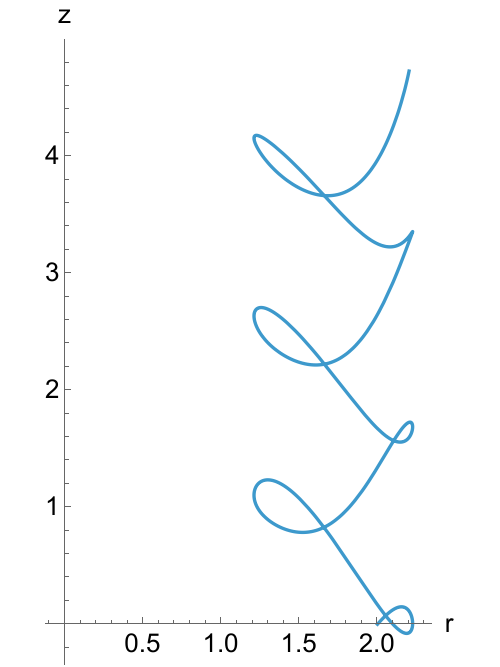}}
  \caption {Plots of $(t,z(t))$ (left), $(r(t),z(t))$ (right) and  $(\|z_1(t)\|,\|z_2(t)\|, z(t))$  (centre) of a sub-Riemannian geodesic passing through point $(0,-\sqrt{2},\sqrt{2},0,0)$.}
  \label{fig: ring geodesic}
\end{figure}

\end{example}

\par\textbf{Acknowledgment.}
The author expresses their sincere gratitude to the anonymous reviewer for the exceptionally careful reading of the manuscript and for the many valuable comments and suggestions that significantly improved the quality and presentation of the paper.

\par\textbf{Funding.}
This research  is partially supported by the Ministry of  Science, Technological Development and Innovation, Republic of Serbia, through the project 451-03-33/2026-03/200104.

\end{document}